\documentclass[a4paper,10pt]{article}

\usepackage{amssymb,amsmath,tikz,amsthm,xypic}

\usepackage{enumerate}

 \usepackage[left=2.6cm, right=2.6cm, top=2.7cm, bottom=3cm]{geometry}

\newtheorem{lemma}{Lemma}

\theoremstyle{definition}
\newtheorem{example}[lemma]{Example}

\def\N{\mathbb N}
\def\Z{\mathbb Z}
\def\Q{\mathbb Q}
\DeclareMathOperator{\en}{\widetilde{ent}}
\DeclareMathOperator{\ent}{{ent}}

\author{Daniele Toller \and Simone Virili\footnote{The second named author was supported by the Spanish Ministerio de Econom\'ia y Competitividad via a research project (MTM2016-77445-P). He was also supported by the Fundaci\'on S\'eneca of Murcia (19880/GERM/15) with a part of FEDER funds.
}}

\title{ERRATA CORRIGE\\Intrinsic algebraic entropy}

\newcommand{\keywords}[1]{\noindent\textbf{Key words and phrases:} #1}
\newcommand{\MSC}[1]{\noindent\textbf{2010 Mathematics Subject Classification.\,} #1}

\date{}

\begin{document}

\maketitle

\abstract{
The notion of intrinsic algebraic entropy of an endomorphism of a given  Abelian group has been recently introduced in [D.\,Dikranjan, A.\,Giordano Bruno, L.\,Salce, S.\,Virili, {\em Intrinsic algebraic entropy}, J.\,Pure Appl.\,Algebra 219 (2015) 2933--2961]. In this short note we provide a correct argument to prove one of the basic properties of the intrinsic algebraic entropy: the Logarithmic Law. In fact, this property was correctly stated in [op.\,cit.] but, as we will show with an explicit counterexample, the original proof contains a flaw. 
}


\MSC{Primary: 20K30; 20K27; Secondary: 20K15; 22B05; 16D10; 37A35; 11R06.}

\section*{Introduction}
The notion of intrinsic algebraic entropy $\en(\phi)$ of an endomorphism $\phi\colon G\to G$ of a given Abelian group $G$ has been recently introduced in \cite{DGSV}, as a natural extension of the algebraic entropy $\ent(\phi)$, studied in \cite{DGSZ}. We refer to \cite{DGSV,DGSZ} for more details on the history of these and related invariants.

The new invariant $\en$ has the advantage that, contrarily to $\ent$, it is non-trivial on torsion-free groups and, furthermore, it satisfies all the desirable properties one should expect from a well-behaved entropy function: it is monotone under taking quotients and restrictions and, in fact, it is additive (see \cite[Thm.\,5.9]{DGSV}), it is upper-continuous (see \cite[Lem.\,3.14]{DGSV}), it satisfies the Logarithmic Law, it takes the ``correct values" on Bernoulli Shifts (see \cite[Ex.\,3.7]{DGSV}), it can be computed easily on endomorphisms of (subgroups of) $\Q^n$ via the so-called Intrinsic Yuzvinski Formula (see \cite[Thm.\,4.2]{DGSV}). Let us remark that the proofs of two of the main properties of the intrinsic entropy, namely its additivity and the Intrinsic Yuzvinski Formula, have been considerably simplified in \cite{GV2} and \cite{SV2}.

\medskip
In this short note we concentrate on the Logarithmic Law for the intrinsic algebraic entropy (see below for the precise definitions and statements). This logarithmic behavior is one of the typical properties for an entropy function, and this is probably the reason for which nobody noticed that its original proof was based on the false lemma \cite[Lem.\,3.11(b)]{DGSV}. 
The situation changed when the first named author, while trying to adapt some of the proofs in \cite{DGSV} to a more general context, noticed the mistake and communicated it to the authors of \cite{DGSV}. In this note we give a simple and explicit counterexample to a wrong statement in \cite[Lem.\,3.11(b)]{DGSV} and we provide a new argument to prove the Logarithmic Law for $\en$. 

\subsection*{Acknowledgement}
It is a pleasure for us to thank Dikran Dikranjan, Anna Giordano Bruno and Luigi Salce for the useful discussions and comments on the matters of this note.

\subsection*{Notations and conventions}
We denote by $\N$ the natural numbers, by $\N_+=\N\setminus\{0\}$ the positive integers, and by
$\log(n)$ the natural logarithm of a positive integer $n$. Throughout this note, $G$ is an  Abelian group, and $\phi\colon G \to G$ is a group endomorphism. 

\subsection*{Basic definitions about the intrinsic algebraic entropy}
For a subgroup $H$ of $G$ and $n\in\N_+$, the {\bf $n$-th partial $\phi$-trajectory of $H$} is the subgroup
\begin{equation*}
T_n(\phi,H) := H + \phi(H)+ \phi^2(H)+ \ldots + \phi^{n-1}(H),
\end{equation*}
while the {\bf $\phi$-trajectory of $H$} is
\begin{equation*}
T(\phi,H) := \bigcup_{n\in\N_+} T_n(\phi,H).
\end{equation*}
Note that $T(\phi,H)$ is the smallest $\phi$-invariant subgroup of $G$ containing $H$.

\medskip
A subgroup $H$ of $G$ is {\bf $\phi$-inert} if the quotient group $( H+\phi(H) ) /H$ is finite \cite[Def.\,1.1]{DGSV}. According to  \cite[Lem.\,3.2]{DGSV}, for such a subgroup, the limit
\begin{equation*}
\en(\phi,H) := \lim_{n \to \infty} \frac{ \log |T_n(\phi,H)/H | }{n}
\end{equation*}
exists, it is finite, and it is called the {\bf intrinsic entropy of $\phi$ with respect to $H$}. Finally, the {\bf intrinsic entropy of $\phi$} is defined in \cite[Def.\,1.2]{DGSV} as the supremum
\begin{equation*}
\en(\phi) := \sup \{ \en(\phi,H) \mid H \text{ $\phi$-inert subgroup of } G \}.
\end{equation*}

\subsection*{The Logarithmic Law}
The ``Logarithmic Law'' for the intrinsic entropy is the following result:

\medskip\noindent
{\bf Logarithmic Law} {\rm(\cite[Lem.\,3.12]{DGSV})}{\bf.}
 {\em 
 If $\phi\colon G \to G$ is an endomorphism of an  Abelian group $G$, and $k \in \N_+$, then 
$
\en(\phi^k) = k \cdot \en(\phi).
$ 
}

\medskip
The Logarithmic Law was stated as Lemma 3.12 in \cite{DGSV}. The original argument for the proof of this property in \cite{DGSV} strongly depends on a second, unfortunately wrong, lemma (\cite[Lem.\,3.11]{DGSV}) which we state below:
 
\begin{lemma}[{\rm \cite[Lem.\,3.11]{DGSV}}]\label{lemma 3.11}
Let $\phi\colon G \to G$ be an endomorphism of an  Abelian group $G$, and $H$ be a subgroup of $G$. 
\begin{enumerate}[(a)]
\item If $H$ is $\phi$-inert, and $H' = T_k(\phi,H)$ for some $k \in \N_+$, then $H'$ is $\phi$-inert and 
\[\en(\phi,H) = \en(\phi, H').\]
\item If $H$ is $\phi^k$-inert for some $k \in \N_+$, and $H' = T_k(\phi,H)$, then $H'$ is $\phi$-inert (so, in particular, $\phi^k$-inert) and $\en(\phi^k,H) = \en(\phi^k, H')$.
\end{enumerate}
\end{lemma}

The last equality stated in Lemma~\ref{lemma 3.11}(b) does not hold true in general, as we show in Example~\ref{example}. As we already mentioned above, Lemma~\ref{lemma 3.11} is used in the proof of \cite[Lem.\,3.12]{DGSV}, in particular,  part (b) of Lemma~\ref{lemma 3.11} is essential in proving the following inequality:
\begin{equation}\label{1/2lg}
\en(\phi^k) \leq k \cdot \en(\phi).
\end{equation}
On the other hand, the second part of the proof of \cite[Lem.\,3.12]{DGSV}, establishing that
\begin{equation}\label{half logarithmic law}
\en(\phi^k) \geq k \cdot \en(\phi),
\end{equation}
is independent of Lemma \ref{lemma 3.11}(b), and it is therefore correct.

\section{Proof of the Logarithmic Law}

We now  give our proof of the Logarithmic Law; for this we will use some deep properties of the intrinsic entropy proved in \cite{DGSV} (more precisely, we will use \cite[Prop.\,5.6, Lem.\,2.7, Prop.\,3.16(b), Lem.\,3.14]{DGSV}), whose proofs do not rely on Lemma \ref{lemma 3.11}(b).

\begin{proof}[Proof of the Logarithmic Law.]
If $\en(\phi) = \infty$, then (\ref{half logarithmic law}) is sufficient to conclude, so let us assume  that $\en(\phi) < \infty$. 

We first prove the Logarithmic Law under the additional assumption that $G = T(\phi, F)$, for a finitely generated subgroup $F\leq G$. Indeed, by \cite[Prop.\,5.6]{DGSV} there exists $m \in \N_+$ such that $T_m(\phi,F)$ is $\phi$-inert and, by Lemma \ref{lemma 3.11}(a),  also the subgroup $H:=T_{m+k-1}(\phi,F)=T_k(\phi,T_m(\phi,F))$ is  $\phi$-inert. Moreover, 
%
%
%
by \cite[Lem.\,2.7(b)]{DGSV}, $H$ is also $\phi^k$-inert. Furthermore, $H$ is a finitely generated subgroup of $G$ such that $G=T(\phi,H)=T(\phi^k,H)$. Then, \cite[Prop.\,3.16(b)]{DGSV} applies respectively to $\phi$ and $\phi^k$ to give the following equalities:
\[
\en(\phi)=\en(\phi,H), \text{ and } \en(\phi^k)=\en(\phi^k,H).
\] 
To conclude this part of the proof, we now verify that $\en(\phi^k,H) = k \cdot \en(\phi,H)$. 
Indeed, for every $n \in \N_+$, $T_n(\phi^k, H) = T_{kn-k+1}(\phi,H)$,
and so 
\begin{align*}
\en(\phi^k, H )&= \lim_{n\to\infty} \frac{ \log| T_n(\phi^k, H) / H | }{n}\\
&= \lim_{n\to\infty} \frac{ \log|  T_{kn-k+1}(\phi,H) /H | }{kn-k+1}\cdot \frac{kn-k+1}{n} 
= k \cdot \en(\phi, H ).
\end{align*}

In the general case, we use the fact that $G$ is the direct limit of the family  of $\phi$-invariant subgroups $\{ T(\phi,F) \mid F \leq G, \ F \text{ finitely generated} \}$. 
Then, \cite[Lem.\,3.14]{DGSV} implies that
\begin{equation}\label{2}
\begin{split}
\en(\phi) &= \sup \{ \en( \phi\restriction_{ T(\phi,F) } )  \mid F \leq G, \ F \text{ finitely generated} \}, \text{ and} \\	
\en(\phi^k) &=  \sup \{ \en( \phi^k\restriction_{ T(\phi,F) } )  \mid F \leq G, \ F \text{ finitely generated} \}.		\end{split}
\end{equation}
For any finitely generated subgroup $F\leq G$, the first part of the proof gives that
$\en( \phi^k\restriction_{ T(\phi,F) } ) = k \cdot \en( \phi\restriction_{ T(\phi,F) } )$,
so (\ref{2}) yields $\en(\phi^k) = k \cdot \en(\phi)$.
\end{proof}

Now we provide a counterexample to the last assertion of Lemma \ref{lemma 3.11}(b).
\begin{example}\label{example}
Let $\Z(2)=\{0,1\}$ be the group of integers modulo $2$, and 
\[
\xymatrix{G := \bigoplus_{i \in \N} \Z(2) =\left\{ (x_i)_{i \in \N} \in \prod_{i \in \N} \Z(2) \mid x_i \neq 0 \text{ for finitely many } i\in \N \right\}.}
\] 
Consider the right Bernoulli shift $\beta\colon G \to G$, mapping $(x_0, x_1, x_2, \ldots ) \mapsto (0, x_0, x_1, x_2, \ldots)$. 
Let $H := \{ (x_i)_{i \in \N} \in G \mid x_i = 0 \text{ for } i \neq 0 \}$, and $H' := T_2(\beta, H)$, so that 
\[
H'= \{ (x_i)_{i \in \N} \in G \mid x_i = 0 \text{ for } i \neq 0,\, 1\}.
\]
 Both $H$ and $H'$ are finite, hence $\beta$-inert and $\beta^2$-inert subgroups of $G$. Now, for a positive integer $n$, 
\[
T_n(\beta^2, H') = \{ (x_i)_{i \in \N} \in G \mid x_i = 0 \text{ for } i \geq 2n \},
\]
so $|T_n(\beta^2, H')/H'|=2^{2n-2}$ and $\en(\beta^2, H') = 2\cdot \log(2)$.
On the other hand, 
\[
T_n(\beta^2, H) = \{ (x_i)_{i \in \N} \in G \mid x_i = 0 \text{ for either } i>2n-2, \text{ or } i \text{ even} \},
\]
so $|T_n(\beta^2, H)/H|=2^{n-1}$ and $\en(\beta^2, H) = \log(2)$.
In particular, $\en(\beta^2, H) \neq \en(\beta^2,  H') $.
\end{example}

\bigskip
--------------------------------

\smallskip
Daniele Toller -- {\tt daniele.toller@uniud.it}\\
Dipartimento di Matematica e Informatica, Universit\`{a} di Udine, Via delle Scienze  206, 33100 Udine, Italy.

\smallskip
Simone Virili -- {\tt virili.simone@gmail.com}\ \  or\ \  {\tt s.virili@um.es}\\
Departamento de Matem\'aticas, Universidad de Murcia, Aptdo.\,4021, 30100 Espinardo, Murcia, Spain.

\end{document}